%\input /texwork/mat/mat3/neumac.tex
%\input /texwork/mat/mat3/neumac
%\catcode`\^^Z=9
%\catcode`\^^M=10
%\input c:/texwork/mat/mat3/german
%
\output={\if N\header\headline={\hfill}\fi
\plainoutput\global\let\header=Y}
\magnification\magstep1
\tolerance = 500
\hsize=14.4true cm
\vsize=22.5true cm
\parindent=6true mm\overfullrule=2pt
\nopagenumbers
\newcount\kapnum \kapnum=0
\newcount\parnum \parnum=0
\newcount\procnum \procnum=0
\newcount\nicknum \nicknum=1
\font\ninett=cmtt9

\font\ninebf=cmbx9

\font\sixbf=cmbx6
\font\ninesl=cmsl9

\font\nineit=cmti9

\font\ninerm=cmr9

\font\sixrm=cmr6
\font\ninei=cmmi9
\font\eighti=cmmi8
\font\sixi=cmmi6
\skewchar\ninei='177 \skewchar\eighti='177 \skewchar\sixi='177
\font\ninesy=cmsy9
\font\eightsy=cmsy8
\font\sixsy=cmsy6
\skewchar\ninesy='60 \skewchar\eightsy='60 \skewchar\sixsy='60
\font\titelfont=cmr10 scaled 1440
\font\paragratit=cmbx10 scaled 1200

\font\name=cmcsc10
\font\emph=cmbxti10

\font\tenmsbm=msbm10
%\font\ninemsbm=msbm9
\font\sevenmsbm=msbm7
%\font\sixmsbm=msbm6
%\font\fivemsbm=msbm5
%\textfont\extsym=\tenmsbm
%\scriptfont\extsym=\sevenmsbm
%\scriptscriptfont\extsym=\fivemsbm
%

%
\font\got=eufm10
\font\Got=eufm7
\font\teneufm=eufm10
\font\seveneufm=eufm7
\font\fiveeufm=eufm5
\newfam\eufmfam
\textfont\eufmfam=\teneufm
\scriptfont\eufmfam=\seveneufm
\scriptscriptfont\eufmfam=\fiveeufm

\font\tenmsam=msam10
\font\sevenmsam=msam7
\font\fivemsam=msam5
\newfam\msamfam
\textfont\msamfam=\tenmsam
\scriptfont\msamfam=\sevenmsam
\scriptscriptfont\msamfam=\fivemsam
\font\tenmsbm=msbm10
\font\sevenmsbm=msbm7
\font\fivemsbm=msbm5
\newfam\msbmfam
\textfont\msbmfam=\tenmsbm
\scriptfont\msbmfam=\sevenmsbm
\scriptscriptfont\msbmfam=\fivemsbm
\def\Bbb#1{{\fam\msbmfam\relax#1}}
\def\cz{{\kern0.4pt\Bbb C\kern0.7pt}
}
\def\ez{{\kern0.4pt\Bbb E\kern0.7pt}
}
\def\fz{{\kern0.4pt\Bbb F\kern0.3pt}}
\def\gz{{\kern0.4pt\Bbb Z\kern0.7pt}}
\def\hz{{\kern0.4pt\Bbb H\kern0.7pt}
}
\def\kz{{\kern0.4pt\Bbb K\kern0.7pt}
}
\def\nz{{\kern0.4pt\Bbb N\kern0.7pt}
}
\def\oz{{\kern0.4pt\Bbb O\kern0.7pt}
}
\def\rz{{\kern0.4pt\Bbb R\kern0.7pt}
}
\def\sz{{\kern0.4pt\Bbb S\kern0.7pt}
}
\def\pz{{\kern0.4pt\Bbb P\kern0.7pt}
}
\def\qz{{\kern0.4pt\Bbb Q\kern0.7pt}
}
\newskip\ttglue
\def\ninepoint{\def\rm{\fam0\ninerm}%
\textfont0=\ninerm \scriptfont0=\sixrm \scriptscriptfont0=\fiverm
\textfont1=\ninei \scriptfont1=\sixi \scriptscriptfont1=\fivei
\textfont2=\ninesy \scriptfont2=\sixsy \scriptscriptfont2=\fivesy
\textfont3=\tenex \scriptfont3=\tenex \scriptscriptfont3=\tenex
\def\it{\fam\itfam\nineit}%
\textfont\itfam=\nineit
\def\sl{\fam\slfam\ninesl}%
\textfont\slfam=\ninesl
\def\bf{\fam\bffam\ninebf}%
\textfont\bffam=\ninebf \scriptfont\bffam=\sixbf
\scriptscriptfont\bffam=\fivebf
\def\tt{\fam\ttfam\ninett}%
\textfont\ttfam=\ninett
\tt \ttglue=.5em plus.25em minus.15em
\normalbaselineskip=11pt
\font\name=cmcsc9
\let\sc=\sevenrm
\let\big=\ninebig
\setbox\strutbox=\hbox{\vrule height8pt depth3pt width0pt}%
\normalbaselines\rm
\def\sl{\it}}

\headline={\ifodd\pageno\rightheadline\else\leftheadline\fi}
\def\rightheadline{\ninepoint Paragraphen"uberschrift\hfill\folio}
\def\leftheadline{\ninepoint\folio\hfill Chapter"uberschrift}
\let\header=Y
\def\titel#1{\need 9cm \vskip 2truecm
\parnum=0\global\advance \kapnum by 1
{\baselineskip=16pt\lineskip=16pt\rightskip0pt
plus4em\spaceskip.3333em\xspaceskip.5em\pretolerance=10000\noindent
\titelfont Chapter \uppercase\expandafter{\romannumeral\kapnum}.
#1\vskip2true cm}\def\leftheadline{\ninepoint
\folio\hfill Chapter \uppercase\expandafter{\romannumeral\kapnum}.
#1}\let\header=N
}
\def\Titel#1{\need 9cm \vskip 2truecm
\global\advance \kapnum by 1
{\baselineskip=16pt\lineskip=16pt\rightskip0pt
plus4em\spaceskip.3333em\xspaceskip.5em\pretolerance=10000\noindent
\titelfont\uppercase\expandafter{\romannumeral\kapnum}.
#1\vskip2true cm}\def\leftheadline{\ninepoint
\folio\hfill\uppercase\expandafter{\romannumeral\kapnum}.
#1}\let\header=N
}
\def\need#1cm {\par\dimen0=\pagetotal\ifdim\dimen0<\vsize
\global\advance\dimen0by#1 true cm
\ifdim\dimen0>\vsize\vfil\eject\noindent\fi\fi}
\def\neupara#1{\par\penalty-2000
\procnum=0\global\advance\parnum by 1
\vskip1cm\noindent{\paragratit \the\parnum. #1}%
\def\rightheadline{\ninepoint\S\the\parnum.\ #1\hfill \folio}%
\vskip 8mm\noindent}
\def\Proclaim #1 #2\finishproclaim {\bigbreak\noindent
{\bf#1\unskip{}. }{\it#2}\medbreak\noindent}
%
%\llap{#3\ }
\gdef\proclaim #1 #2 #3\finishproclaim {\bigbreak\noindent%
\global\advance\procnum by 1
{%
{\relax\ifodd \nicknum
\hbox to 0pt{\vrule depth 0pt height0pt width\hsize
\quad \ninett#3\hss}\else {}\fi}%
\bf\the\parnum.\the\procnum\ #1\unskip{}. }
{\it#2}
\immediate\write\num{\string\def
\expandafter\string\csname#3\endcsname
{\the\parnum.\the\procnum}}
\medbreak\noindent}
\newcount\stunde \newcount\minute \newcount\hilfsvar
\def\uhrzeit{
\stunde=\the\time \divide \stunde by 60
\minute=\the\time
\hilfsvar=\stunde \multiply \hilfsvar by 60
\advance \minute by -\hilfsvar
\ifnum\the\stunde<10
\ifnum\the\minute<10
0\the\stunde:0\the\minute~Uhr
\else
0\the\stunde:\the\minute~Uhr
\fi
\else
\ifnum\the\minute<10
\the\stunde:0\the\minute~Uhr
\else
\the\stunde:\the\minute~Uhr
\fi
\fi
}
\def\calA{{\cal A}} 
\def\calC{{\cal C}} 
 
 \def\calH{{\cal H}}

 \def\gotA{\hbox{\got A}}

\def\goth{\hbox{\got h}}

\def\gotr{\hbox{\got r}}

\def\gotw{\hbox{\got w}} 
 
\def\goty{\hbox{\got y}} 
\def\gotz{\hbox{\got z}} 

 \def\Gotz{\hbox{\Got z}}
\def\schreib#1{\hbox{#1}}

\def\mod{\mathop{\rm mod}\nolimits}
\def\O{{\rm O}}

\def\re{\mathop{\rm Re}\nolimits}
\def\Re{\re}

\def\SO{\mathop{\rm SO}\nolimits}

\def\Sp{\mathop{\rm Sp}\nolimits}

\def\Spin{\mathop{\rm Spin}\nolimits}

\def\boxit#1{
\vbox{\hrule\hbox{\vrule\kern6pt
\vbox{\kern8pt#1\kern8pt}\kern6pt\vrule}\hrule}}
\def\Boxit#1{
\vbox{\hrule\hbox{\vrule\kern2pt
\vbox{\kern2pt#1\kern2pt}\kern2pt\vrule}\hrule}}

\def\smallni{\smallskip\noindent }

\def\bigni{\bigskip\noindent }
\def\Isom{\mathop{\;{\buildrel \sim\over\longrightarrow }\;}}
\def\lo{\longrightarrow}

\def\loma{\longmapsto}

\def\pii{\pi {\rm i}}

\def\square{\hbox{\hbox to 0pt{$\sqcup$\hss}\hbox{$\sqcap$}}}
\def\qed{\ifmmode\square\else{\unskip\nobreak\hfil
\penalty50\hskip3em\null\nobreak\hfil\square
\parfillskip=0pt\finalhyphendemerits=0\endgraf}\fi}
\def\pn{\the\parnum.\the\procnum}
\def\downmapsto{{\buildrel
{\vbox{\hbox{\hskip.2pt$\scriptstyle-$}}}
\over{\raise7pt\vbox{\vskip-4pt\hbox{$\textstyle\downarrow$}}}}}
\input eachttheta.num \newwrite\num\openout\num=eachttheta.num
\nicknum=0      %bewirkt, dass keine Nicknamen am Rand
\def\RAND#1{\hbox to 0mm{\hss\vtop to 0pt{%
\raggedright\ninepoint\parindent=0pt%
\baselineskip=1pt\hsize=2cm #1\vss}}}
\noindent\centerline{\titelfont Octavic theta series}%
\def\tr{{\rm tr}} 
  
\vskip 1.5cm
\leftline{\it \hbox to 6cm{Eberhard Freitag\hss}
Riccardo Salvati
Manni  }
\leftline {\it  \hbox to 6cm{Mathematisches Institut\hss}
Dipartimento di Matematica, }
\leftline {\it  \hbox to 6cm{Im Neuenheimer Feld 288\hss}
Piazzale Aldo Moro, 2}
\leftline {\it  \hbox to 6cm{D69120 Heidelberg\hss}
I-00185 Roma, Italy. }
\leftline {\tt \hbox to 6cm{freitag@mathi.uni-heidelberg.de\hss}
salvati@mat.uniroma1.it}
\def\text{\hbox}

\def\leftheadline{\ninepoint\folio\hfill
Octavic modular forms}%
\def\rightheadline{\ninepoint Introduction\hfill \folio}%
\headline={\ifodd\pageno\rightheadline\else\leftheadline\fi}%
\bigni
\centerline{\paragratit \rm  2014}%
\let\header=N%

%FIRST DRAFT
\vskip5mm\noindent
{\paragratit Introduction}%
\bigni
In the paper [FS] we considered the even unimodular lattice $L=\Pi_{2, 10}$ of signature $(2,10)$.
It can be realized as direct sum of the negative of the lattice $E_8$ and two hyperbolic
planes. We considered the orthogonal group $\O(L)$. A certain subgroup of index two
$\O^+(L)$
acts biholomorphically on a ten dimensional tube domain $\calH_{10}$. We denote the variables
with respect to the standard embedding $\calH_{10} \hookrightarrow\cz^{10}=\Pi_{1, 9}\otimes\cz$ which we used
in [FS] by $z$ (compare Sect.~4). We recall one of the main results in [FS].
\Proclaim
{Theorem}
{There exists a non vanishing modular form $f(z)$
of weight $4$ (the singular weight) with respect
to the full modular group $\O^+(L)$. It is uniquely determined up to a
constant factor. The form $f(2z)$ belongs to the principal
congruence subgroup of level two $\O^+(L)[2]$. The $\O^+(L)$-orbit of $f(2z)$ spans a
$715$-dimensional space which is the direct sum of\/ $\cz f(z)$ and a $714$-dimensional
irreducible space.}
\finishproclaim
The existence and uniqueness of $f(z)$ can already be derived from the paper [EK]
of Eie and Krieg.
In [FS] we obtained also the following result.
The $715$-dimensional space defines an everywhere regular birational embedding
of the associated modular variety into $\pz^{714}$. The image is contained in a certain system
of quadrics. If one is very optimistic, one may conjecture that the ring of modular forms
of weight divisible by 4 is generated by the 715-dimensional space and that the quadratic
relations are  the defining ones.
\smallskip
There is a close relation with the Borcherds $\Phi$-function [Bo] which is a 
modular form of weight 4 with respect to the orthogonal group of the Enriques lattice
$\sqrt 2M$ where   
$$M=U\oplus \sqrt2 U\oplus (-\sqrt 2 E_8).$$
In fact, as it is explained in [Bo], the lattice $2L$ and the inverse image 
of any non-zero  isotropic   vector  of $L/2L$ generate a copy of the 
Enriques lattice. Hence each non-zero isotropic vector $\alpha$ of $L/2L$
leads to a realization $f_\alpha$ of the $\Phi$-function inside our 715-dimensional
space. Their divisors are Heegner divisors $H_{\alpha}(-2)$ which belong to vectors of norm
$-2$. They  have an infinite  product expansion.   
In the paper [FS] we proved that the $2079$ forms $f_\alpha$
generate the 714-dimensional space. 
\smallskip
Recently, in [KMY], it has been proved that the restriction of  $\Phi$
to a subdomain isomorphic to the Siegel space $\hz_{2}$ of 
degree 2 can be expressed as a theta series.
As a by-product  they    showed that  the  eighth power of  
any  even theta constant in  genus  two
is expressed as an infinite product of Borcherds type.
\smallskip
So we have been asked  by one  of the  authors  whether $\Phi$ is related to theta series.
In this paper we give an affirmative answer. We will construct a modular embedding
of $\calH_{10}$ into the Siegel half plane $\hz_{16}$ of degree 16. This means that every
substitution of $\O^+(L)$ extends to a Siegel modular substitution. Even more, we will show that
it extends to a substitution of the theta group $\Gamma_{16,\vartheta}$ which is the group
$\Gamma_{16}[1,2]$ in Igusa's notation. As a consequence, $f(z)$ can be constructed as
the restriction of the simplest among all theta series.
$$\vartheta(Z)=\sum_{g\in\gz^n}e^{\pii Z[g]} \qquad (n=16).$$
This modular embedding has  also the property that the transformations in $\O^+(L)[2]$
extend to transformations in  $\Gamma_{16}[2,4]$. The most natural modular forms
on $\Gamma_n[2,4]$ are the theta series of second kind
$$\sum_{g\in\gz^n}e^{2\pii Z[g+a/2]},\quad a\in\gz^n.$$
We will see that their restrictions span the 715-dimensional space.
\smallskip
We have seen already in [FS] that the forms of weight 4 are determined by their values at
the zero dimensional cusp classes. Since one can compute these values in both pictures,
one can make the identification of the orthogonal modular forms in [FS] and the restrictions
of the thetas of second kind explicit. 
\smallskip
The construction of the modular embedding rests on some results about the Clifford algebra
of the lattice $-E_8$ which we will derive by means of the octavic multiplication.
Similar constructions can be found in [FH].
\neupara{The Clifford algebra}%
Let $(V,q)$ be a quadratic space of positive
dimension over a field $K$. We want to include characteristic two, so we recall
that $q$ means a map $q:V\to K$ with the properties
\vskip1mm
\item{a)} $q(ta)=t^2q(a)$ for $t\in K$,
\item{b)}$(a,b)=q(a+b)-q(a)-q(b)$ is bilinear.
\smallni
The Clifford algebra $\calC(V)$ is
an associative algebra with unit which contains $V$ as sub-vector space
and which
is generated by $V$ as algebra. The defining relations are
$$ab+ba=(a,b)\qquad(a,b\in V).$$
Of course $K$ is embeddedd into $\calC(V)$ by $t\mapsto t\,1_{\calC(V)}$.
This defines an embedding
$$K\oplus V\lo \calC(V).$$
The main involution of $\calC(V)$ is denoted by $a\mapsto a'$. This is an
involutive antiisomorphism which acts on $V$ as the negative of the identity,
$$(a+b)'=a'+b',\ (ab)'=b'a',\quad a'=-a\schreib{ for } a\in V.$$
The even part of $\calC(V)$ is the subalgebra $\calC^+(V)$ generated by the
two-products $ab$, $a,b\in V$. It is invariant under the main involution.
We are interested in the case $K=\rz$, when the dimension of $V$ is eight
and   $q$ is negative definite. It is known in this case
that there exists an isomorphism of involutive algebras
$$\calC^+(\rz^8)\cong M_8(\rz)\times M_8(\rz)\qquad(q(x)=-x_1^2-\cdots-x_8^2).$$
Here the involution on $M_8(\rz)\times M_8(\rz)$ is defined by
$(A,B)\mapsto (A',B')$, where the dash now means matrix transposition.
As a consequence, there exists an homomorphism of involutive algebras
$$\calC^+(\rz^8)\lo M_8(\rz)\qquad (q(x)=-x_1^2-\cdots-x_8^2).$$
It is important for us to get an explicit homomorphism, which preserves also
a certain integral structure. This comes into the game if we consider
a lattice $(L,q)$ and $V=L\otimes_\gz\rz$. Then we can consider the
$\gz$-algebra $\calC(L)\subset\calC(V)$ generated by $L$ and also the
$\gz$-algebra $\calC^+(L)$ generated by all $ab$, $a,b\in L$.
Both algebras are invariant under the main involution.
We are in particular interested in the lattice $E_8$
and its negative definite version. Sometimes we write simply
$-E_8$ for it. We mentioned that there exists an involutive homomorphism
$\calC(-E_8\otimes_\gz\rz)\to M_8(\rz)$. We want
that it induces a (surjective) homomorphism
$$\calC(-E_8)\lo M_8(\gz).$$
We will use the octavic multiplication for this purpose.
\neupara{Octaves}%
We consider $\oz=\rz^8$ with the standard basis
$$e_0=(1,0,0,0,0,0,0,0),\dots, e_7=(0,0,0,0,0,0,0,1).$$
The octavic product  is the bilinear
map
$$\oz\times\oz\to\oz,\quad (a,b)\loma a*b,$$
defined by the multiplication table
$$\matrix{
e_0&e_1&e_2&e_3&e_4&e_5&e_6&e_7\cr
e_1&-e_0&e_4&e_7&-e_2&e_6&-e_5&-e_3\cr
e_2&-e_4&-e_0&e_5&e_1&-e_3&e_7&-e_6\cr
e_3&-e_7&-e_5&-e_0&e_6&e_2&-e_4&e_1\cr
e_4&e_2&-e_1&-e_6&-e_0&e_7&e_3&-e_5\cr
e_5&-e_6&e_3&-e_2&-e_7&-e_0&e_1&e_4\cr
e_6&e_5&-e_7&e_4&-e_3&-e_1&-e_0&e_2\cr
e_7&e_3&e_6&-e_1&e_5&-e_4&-e_2&-e_0\cr}$$
It is known that this product has no zero divisors.
We will call the elements of $\oz$ from now on ``octaves''.
An octave is  integral, if it is in the $\gz$-module
$$\oz(\gz):=\gz f_0+\cdots+\gz f_7,$$
where
$$f_0=e_0,\quad  f_1=e_1,\quad  f_2=e_2,\quad f_3=e_3$$
and
$$\eqalign{f_4={e_1+e_2+e_3-e_4\over 2},\quad &f_5={-e_0-e_1-e_4+e_5\over 2},\cr
f_6={-e_0+e_1-e_2+e_6\over 2},\quad &f_7={-e_0+e_2+e_4+e_7\over 2}.\cr}$$
This is a subring (closed under octavic multiplication). We also extend the
octavic multiplication $\cz$-linearly
to $\oz(\cz):=\cz^8$.
\smallskip
The conjugate of an octave $x\in\oz$ is defined by
$$\overline{x_0e_0+x_1e_1+\cdots+x_7e_7}=x_0e_0-x_1e_1+\cdots-x_7e_7.$$
One has
$$\overline{x*y}=\bar y* \bar{x}.$$
Because the element $e_0$ is the unit element, we embed
$\rz$ into $\oz$ by sending $t$ do $te_0$. Sometimes we identify $t$ with its
image $te_0$.
The norm is defined by
$$N:\rz^8\lo\rz,\quad N(x)=\bar x x=x_0^2+\cdots+x_7^2.$$
One has $N(xy)=N(x)N(y)$.
The norm of an integral octave is an integer.
We consider also the trace
$$\tr:\rz^8\lo\rz,\quad \tr(x)=x+\bar x=2x_0.$$
Sometimes $\Re(x):=x_0=\tr(x)/2$ is called the {\it real part\/} of $x$.
The traces of integral octaves are integral and the trace can be extended
$\cz$-linearly to $\oz(\cz)$.
\smallskip
The norm is a quadratic form with associated bilinear form
$\tr(\bar x * y)$. Hence $(\oz(\gz),N)$ can be considered as an
even lattice. Actually it is a copy of the $E_8$-lattice.
\smallskip
Since the octavic multiplication is not associative, one has sometimes to be
a little careful. On can check
$$\tr((a*b)*c)=\tr(a*(b*c))$$
and then use the notation $\tr(a*b*c)$ for this expression. One also has
$$\tr(a*b*c)=\tr(b*c*a).$$
We consider now the Clifford algebra of $\rz^8$ equipped with the negative definite
form $-x_1^2-\cdots-x_8^2$. We identify $\rz^8=\oz$
and write $-\oz$ to indicated that we consider the negative definit quadratic
form $-N$. We embed $\oz$ into the even part of the Clifford algebra
$\calC^+(-\oz)$,
$$\oz\lo\calC^+(-\oz),\quad a\loma e_0a.$$
The algebra $\calC^+(-\oz)$ is generated as algebra by the image of this map.
Hence a homomorphism starting from this algebra is determined if we know
its values at elements of the form $e_0a$.
\proclaim
{Proposition}
{For an octave $a$ we define the\/ $8\times8$-matrix
$$P(a):=(\tr(\bar e_i*a*e_k)/2).$$
There is a unique involutive homomorphism of algebras
$$\calC^+(-\oz)\lo M_8(\rz),\quad e_0a\loma P(a).$$
}
ShC%
\finishproclaim
{\it Proof.\/} Let be $P_i=P(e_i)$.
The matrix $P_0$ is the negative unit matrix.
The remaining defining relations are
$P_1^2=\cdots=P_7^2=P_0$ and $P_iP_j=-P_jP_i$ for $1\le i<j\le 7$.
They are easy to check.
Hence we obtain a homomorphismus. It is involutive because the $P_1,\dots,P_7$
are skew-symmetric.
\qed
\smallskip
The integral part $\calC^+(-\oz(\gz))$ of the Clifford algebra is the abelian
group generated by the products (including the empty product)
$$f_{i_1}\cdots f_{i_m},\quad 0\le i_1<\cdots<i_m\le7,\quad m\equiv0\mod 2.$$
From the Clifford relations follows
$(e_0f_i)(e_0f_j)=f_if_j+(e_0,f_i)e_0f_j$.
Hence another basis is given by the
$$(f_0f_{j_1})\cdots (f_0f_{j_m}),\quad 1\le j_1<\cdots<j_m\le 7.$$
Hence we obtain
\proclaim
{Corollary}
{The image
of $\calC^+(-\oz(\gz))$ in $M_8(\rz)$ is the  abelian group generated by
the matrices (the empty product included)
$$P(f_{j_1})\cdots P(f_{j_m}),\quad 1\le j_1<\cdots<j_m\le 7.$$}
IgO%
\finishproclaim
\neupara{A spin group}%
We consider now the 12-dimensional vector space $V$ with the basis
$$h_1,h_2,h_3,h_4,\ e_0,\dots,e_7.$$
We equip it with the bilinear form
defined by the
Gram matrix
$$\pmatrix{0&1&&&&&\cr1&0&&&&&\cr&&0&1&&&\cr&&1&0&&&\cr
&&&&-2&&\cr&&&&&\ddots&\cr&&&&&&-2\cr}.$$
By means of the subspaces
$$H_1(\rz):=\rz h_1+\rz h_2,\ H_2(\rz):=\rz h_3+\rz h_4,\
V_0:=\rz e_0+\cdots+\rz e_{7}$$
we obtain an orthogonal decomposition
$$V=H_1(\rz)\oplus H_2(\rz)\oplus V_0.$$
The following description of the structure of the Clifford algebra is
taken from [FH].
\smallskip
Let $\gotA$ by an involutive algebra, i.e.\ an unital associative
$\rz$-algebra which is equipped with an involution
$a\mapsto a'$. We extend this involution to matrices
with entries from $\gotA$ by the
formula
$$(m_{ij})^*:=(m'_{ji}).$$
\proclaim
{Lemma}
{There exists a commutative diagram of algebras
$$\matrix{\calC(V)&\Isom&M_4(\calC(V_0))\cr
\cup&&\cup\cr
\calC^+(V)&\Isom&M_4(\calC^+(V_0)).}$$
The main involution of  $\calC(V)$
corresponds to the involution
$$\pmatrix{a&b\cr c&d}\longmapsto
\pmatrix{d^*&-b^*\cr -c^*&a^*}.$$
(The blocks $a,b,c,d$ are $2\times2$-matrices).
The diagram is defined by the assignments
$$\leqalignno{
h_1&\longmapsto e_0\pmatrix{0&1&0&0\cr0&0&0&0\cr0&0&0&0\cr0&0&1&0},\cr
h_2&\longmapsto e_0\pmatrix{0&0&0&0\cr-1&0&0&0\cr0&0&0&-1\cr0&0&0&0},\cr
h_3&\longmapsto e_0\pmatrix{0&0&0&1\cr0&0&-1&0\cr0&0&0&0\cr0&0&0&0},\cr
h_4&\longmapsto e_0\pmatrix{0&0&0&0\cr0&0&0&0\cr0&1&0&0\cr-1&0&0&0},\cr
e_0&\longmapsto e_0\pmatrix{-1&0&0&0\cr0&1&0&0\cr0&0&-1&0\cr0&0&0&1},\cr
e_i&\longmapsto e_i\pmatrix{1&0&0&0\cr0&1&0&0\cr0&0&1&0\cr0&0&0&1},
\quad 1\le i\le 7.\cr}$$
}
spI
\finishproclaim
{\it Proof of lemma \spI}.\
The images of the basis elements satisfy the defining
relations of the Clifford algebra.
By the universal property of the Clifford algebra
this extends to a
homomorphism $\calC(V)\to M_4(\calC(V_0))$.
One verifies that this homomorphism is surjective.
A dimension argument shows that
it is an isomorphism. The rest is clear.\qed
\smallskip
The {\it spin-group\/} $\Spin(V)$ of a quadratic space $V$
consists of all $g\in\calC^+(V)$
with the properties
$$g'g=1\quad\schreib{and}\quad gVg^{-1}=V.$$
The resulting transformations of $V$
$$x\mapsto gxg^{-1}$$
are orthogonal. This defines a  two to one  homomorphism
$$\Spin(V)\longrightarrow \O(V).$$
The image is the  spinor kenel (i.e. the connected component of the orthogonal  group).
If $L\subset V$ is a lattice we can define the integral spin group
$$\Spin(L):=\Spin(V)\cap\calC^+(L).$$
Using the description of \spI\ of the Clifford algebra of our 12-dimensional
space $V$ we obtain a relation between its spin group and a symplectic group:
\smallskip
The  (Hermitian)
symplectic group (of degree two)
$\Sp(2,\gotA)$
of an involutive algebra $\gotA$
consists of all
$4\times4$-matrices $M$ with entries from $\gotA$, such that
$$M^*IM=I,\qquad I=\pmatrix{0&e\cr-e&0}\qquad
(e=\hbox{\rm unit matrix}).$$
From \spI\ we obtain
\proclaim
{Lemma}
{The restriction of the homomorphism $\spI$ defines an embedding
$$\Spin(V)\lo \Sp(2,\calC^+(V_0)).$$
Using the realization $V_0=-\oz$ and combining with the homomorphism \ShC\
we obtain an injective  homomorphim
$$J_0:\Spin(V)\lo\Sp(16,\rz).$$
Here $\Sp(16,\rz)$ denotes the standard symplectic group of degree $16$
($32\times 32$-matrices).
}
eBS%
\finishproclaim
\neupara{An embedding into the Siegel half plane}%
Again we consider the vector space $V=\rz^{12}=\rz^4\times V_0$ with the basis
$h_1,\dots,h_4$, $e_0,\dots,e_7$ equipped with the quadratic form $q$
of signature $(2,10)$. We denote the corresponding orthogonal half plane
by $\calH_{10}$. It can be defined a the set of triples $(z_1,z_2,\gotz)$, where
$z_1,z_2$ are in the usual upper half plane and where $\gotz\in V_0(\cz)$
such that $y_1y_2+q(\goty)>0$. We recall that $q$ is negative definite on $V_0$.
There is an embedding
$\calH_{10}\to \pz(V(\cz))$ such that the image is one connected
component of the subset  $\tilde\calH_{10}$ defined  
in $V(\cz)$ by $(z,z)=0$ and $(z,\bar z)>0$.
\smallskip
From this we obtain an action of a subgroup of index two $\O^+(V)$ of $\O(V)$.
The image of $\Spin(V)$ is the connected component of $\O(V)$. hence we have
a natural map $\Spin(V)\to\O^+(V)$. We use this map to define an action of
$\Spin(V)$ on $\calH_{10}$. Our next goal is to construct an embedding
of $\calH_{10}$ into the Siegel half plane $\hz_{16}$ of degree 16 which is
compatible with the homomorphism $\Spin(V)\to \Sp(16,\rz)$ and
the standard action of the symplectic group one the Siegel half plane.
For this purpose it is convenient to push $\calH_{10}$ into the even part of the
Clifford algebra, more precisely we consider the map
$$\calH_{10}\lo M_2(\calC^+(V_0(\cz)),\quad (z_1,z_2,\gotz)\loma
\pmatrix{z_1&e_0\gotz\cr(e_0\gotz)'&z_2}.$$
We explain the notations: The real bilinear form $(\cdot,\cdot)$ extends
to a $\cz$-bilinear form on the complexification $V_0(\cz)$ of $V_0$. Its
(complex) Clifford algebra is $\calC(V_0(\cz))$. It can be identified with
$\calC(V_0)\otimes_\rz\cz$. The same is true
for the even part $\calC^+$.
The main involution extends to a $\cz$-linear
involution of $\calC^+(V_0(\cz))$, which is denoted by the same letter.
Now we consider the homomorphism $\calC^+(V_0)\to M_8(\rz)$.
We extend in the natural way to a map
$$M_2(\calC^+(V_0(\cz))\lo M_{16}(\cz).$$
\proclaim
{Proposition}
{The image of\/ $\calH_{10}$ under the maps
$$\calH_{10}\lo M_2(\calC^+(V_0(\cz))\lo M_{16}(\cz)$$
is contained in the Siegel half plane of degree 16. This is an embedding
$j_0:\calH_{10}\to\hz_{16}$
which is compatible with the homomorphism
$J_0:\Spin(V)\to\Sp(16,\rz)$
in the sense that the diagram
$$\matrix{\Spin(V)\times\calH_{10}&\buildrel (J_0,j_0)\over\lo
&\Sp(16,\rz)\times \hz_{16}\cr
\downarrow&&\downarrow\cr
\calH_{10}&\buildrel j_0\over\lo&\hz_{16}\cr}$$
commutes.}
Mee%
\finishproclaim
\neupara{A modular embedding}%
We have to modify the embedding $(J_0,j_0)$
slightly because we want to have that
the integral structures are preserved.
We identify now $V_0=\rz^8$ with the octaves $\oz$ and we consider the lattice
$\oz(\gz)$. In the quadratic  space $V=\rz^4\times(-\oz)$ we consider the lattice
$\gz^4\times(-\oz(\gz))$ which is even and unimodular and of signature $(2,10)$.
Hence we denote it simply by
$$\hbox{\rm II}_{2,10}:=\gz^4\times(-\oz(\gz)).$$
\proclaim
{Lemma}
{The image of $\calC^+(\hbox{\rm II}_{2,10})$
under the homomorphism
\spI\ is
$M_4(\calC^+(-\oz(\gz)))$.}
ICC%
\finishproclaim
We want to modify the homomorphism \ShC\ in such a way that the image of
$\calC^+(\oz(\gz))$ consists of integral matrices. For this
we consider the matrix $F$ whose
rows are the vectors $f_0,\dots,f_7$. Then $S:=2FF'$ is an even unimodular
matrix. We consider the symplectic matrix
$$M=\pmatrix{\sqrt 2F&0&0&0\cr0&\sqrt 2F&0&0\cr
0&0&(\sqrt 2F)'^{-1}&0\cr 0&0&0&(\sqrt2F)'^{-1}}.$$
We combine the map $\calH_{10}\to \hz_{16}$ defined a above with the substitution
$Z\mapsto M(Z)$. Because of
$(\sqrt 2F)P(a)(\sqrt2F)'=(\tr(\bar f_iaf_j))$
we obtain:
\proclaim
{Definition}
{For an octave $a$ we define the matrix
$Q(a)=\bigl(\tr(\bar f_ia f_j)\bigr)$.
This gives a map
$$Q:\oz(\gz)\lo M_8(\gz).$$
We also consider
$$j:\calH_{10}\lo \hz_{16},\quad (z_1,z_2,\gotz)\loma
\pmatrix{z_0S&Q(\gotz)\cr Q(\gotz)'&z_2S}$$
and we consider the homomorphism
$$J:\Spin(V)\lo\Sp(16,\rz),$$
which is the composition of the homomorphism described in \Mee\
and the inner automorphism $N\loma M^{-1}NM$  of\/ $\Sp(16,\rz)$.}
DMe%
\finishproclaim
The advantage of this modified embedding is that it is modular, i.e.~it
preserves the modular groups:
\proclaim
{Proposition}
{The diagram
$$\matrix{\Spin(V)\times\calH_{10}&\buildrel(J,j)\over
\lo &\Sp(16,\rz)\times \hz_{16}\cr
\downarrow&&\downarrow\cr
\calH_{10}&\buildrel j\over\lo&\hz_{16}\cr}$$
commutes.
The homomorphism $J$  maps the group
$\Spin(\hbox{\rm II}_{2,10})$  into the Siegel modular group
$\Sp(16,\gz)$.
}
MeS%
\finishproclaim
{\it Proof.\/}
We have to determine the image of the {\it integral\/} spin-group
$\Spin(\hbox{\rm II}_{2,10})$ in the symplectic group. We have to consider the map
$\calC^+(-\oz(\gz))\to M_8(\rz)$ (see \ShC\ and \IgO) .
We denote its image by $\calA$. We have to consider $4\times4$-blocks
with entries in $\calA$ and conjugate this with $M$. The resulting matrix
consists of blocs of the form
$$FAF^{-1},\quad 2 FAF',\quad
{1\over2}F'^{-1}AF^{-1},\quad F'^{-1}AF'\qquad(A\in\calA).$$
We have to show that they are integral. Because $S=2FF'$ is even unimodular,
they all are integral if one is integral. Hence its remains to remind that
$Q(a)=2FP(a)F'$ is integral for all integral octaves (s.\DMe).\qed
\neupara{The theta group}%
In a first step we investigate the
principal congruence subgroups of level two. 
First of all we recall the  principal congruence subgroup of level $l$ in the symplectic case. 
It  consists of all  integral symplectic matrices of degree $2n$ such that ${A\,B\choose C\,D}$  is congruent  to ${1_n\,0\choose 0\,1_n}\, {\rm  mod}\,l.$ We
denote  it by $\Gamma_{n}[l]$. 
\smallskip
There are two ways to define the principal congruence subgroups of level two
in $\Spin(\hbox{\rm II}_{2,10})$. Firstly one can consider
of $\O(\hbox{\rm II}_{2,10})[2]$, which is the subgroup of
$\O(\hbox{\rm II}_{2,10})$ acting trivial on
$\hbox{\rm II}_{2,10}/2\hbox{\rm II}_{2,10}$. We denote its inverse image in
$\Spin(\hbox{\rm II}_{2,10})$ by $\Spin(\hbox{\rm II}_{2,10})[2]$.
Secondly one can consider the reduction of the Clifford algebra mod 2, which
is nothing else but the Clifford algebra of the quadratic space
$\fz_2^{12}$. The reduction homomorphism
$$\calC(\hbox{\rm II}_{2,10})\lo \calC(\fz_2^{12})=
\calC(\hbox{\rm II}_{2,10})\otimes_\gz\fz_2$$
induces a homomorphism
$$\Spin(\hbox{\rm II}_{2,10})\lo\Spin(\fz_2^{12}).$$
Since we are in characteristic two now, the Spin group is a subgroup
of the orthogonal group. More precisely $\Spin(\fz_2^{12})$ is the simple subgroup
of index two of the orthogonal group $\O(\fz_2^{12})$. We see that the kernel of
the reduction homomorphism is $\Spin(\hbox{\rm II}_{2,10})[2]$. This shows that
every element of this group can be written in the form
$e+2a$, where $a$ is an element of the integral Clifford algebra
$\calC(\hbox{\rm II}_{2,10})$. Now the same proof as in the second part of
\MeS\ shows:
\proclaim
{Remark}
{The image of $\Spin(\hbox{\rm II}_{2,10})[2]$ under the homomorphism $J$ is
contained in the principal congruence subgroup of level two  $\Gamma_{16}[2]$.}
IzC%
\finishproclaim
There is  a rather involved refinement of the second statement of
\MeS:
Recall that the Igusa group $\Gamma_n[l,2l]$
consists of all ${A\,B\choose C\,D}$ in $\Gamma_n[l]$
such that $AB'/l$ and
$CD'/l$ have even diagonal. The group $\Gamma_n[1,2]$ is the so-called theta
group.
We claim now:
\proclaim
{Proposition}
{The image of the integral spin-group $\Spin(\hbox{\rm II}_{2,10})$
under the homomorphism $J$ is contained in the theta-group
$\Gamma_{16}[1,2]$.}
IeT%
\finishproclaim
{\it Proof.\/} It is enough to check the images for a system of elements
whose images in $\Spin(\fz_2^{12})$ generate this group. One can take
the following system of elements from $M_4(\calC^+(-\oz(\gz))$
$$\pmatrix{0&0&1&0\cr0&0&0&1\cr-1&0&0&0\cr0&0&-1&0\cr},\quad
\pmatrix{1&0&h_1&\goth\cr0&0&\goth'&h_2\cr0&0&1&0\cr0&0&0&1\cr},\quad
\pmatrix{1&\goth&0&0\cr0&1&0&0\cr0&0&1&0\cr0&0&-\goth'&1\cr}$$
where $h_1,h_2\in\gz$ and $\goth\in e_0\oz(\gz)$.
It is easy to check that they are contained in $\Spin(\hbox{\rm II}_{2,10})$.
It remains to show
that their
images under $J$ are contained in $\Gamma_{16}[1,2]$. From the definition of
$J$ one sees that this means that the matrix $Q(a)$ (see \DMe)
has even diagonal for $a\in\oz(\gz)$. One has to check this for
$a=f_i$. But for $i>0$ this matrix is skew-symmetric and for $i=0$ one obtains
a Gram-matrix for $E_8$ which is an even lattice.\qed
\smallskip
There is an improvement of \IzC.
\proclaim
{Lemma}
{The image of $\Spin(\hbox{\rm II}_{2,10})[2]$ under the homomorphism $J$ is
contained in the Igusa group
$\Gamma_{16}[2,4]$.}
IzD%
\finishproclaim
{\it Proof.\/}
We consider
$$N=\pmatrix{\sqrt 2 E&0\cr 0&\sqrt 2^{-1} E}\in\Spin(V).$$
It has the effect  N(Z)=2Z for $Z\in\calH_{10}$.
We conjugate an arbitrary element $M\in \Spin(L)[2]$ with
$N$. We claim that the element $N^{-1}MN$ is still in
in $\Spin(\hbox{\rm II}_{2,10})$. For this is sufficient to show
that $N^{-1}MN$ is in the Clifford algebra $\calC(\hbox{\rm II}_{2,10})$
if $M\in 2\calC(\hbox{\rm II}_{2,10})$. This is sufficient to check for  generators
of the $\gz$-algebra $\calC(\hbox{\rm II}_{2,10})$, which is easy.
The replacement $M\mapsto N^{-1}MN$ has the effect
$AB'\mapsto AB'/2$.
Using \IeT\ this leads to $M\in \Spin(L)[2,4]$.\qed
\smallskip
We need the relation between the spin- and the orthogonal group.
\proclaim
{Lemma}
{The image of the natural homomorphism
$$\Spin(\hbox{\rm II}_{2,10})\lo \O(\hbox{\rm II}_{2,10})$$
is the subgroup of index four $\SO^+(\hbox{\rm II}_{2,10})$.}
SOn%
\finishproclaim
It is easy to show that elements of the group $\O^+(\hbox{\rm II}_{2,10})[2]$
have determinant one. So we can reformulate the results about the
spin group as follows:
\proclaim
{Lemma}
{The homomorphism $J$ induces homomorphisms
$$\eqalign{\SO^+(\hbox{\rm II}_{2,10})&\lo\Gamma_{16}[1,2]/{\pm E},\cr
\O^+(\hbox{\rm II}_{2,10})[2]&\lo \Gamma_{16}[2,4]/{\pm E}.\cr}$$}
SIh%
\finishproclaim
\neupara{Theta series}%
We consider the standard theta series on $\hz_{16}$
$$\vartheta_\varphi(Z):=\sum_{g\in\gz^{16}}\phi(g)e^{\pii Z[g]/2},
\quad \varphi: (\gz/2\gz)^{16}\lo\cz.$$
Another way to define them is to use the basis
$$\sum_{g\in\gz^{16}}e^{2\pii Z[g+m/2]},\quad m\in(\gz/2\gz)^{16}.$$
We restrict them to $\calH_{10}$ using the modular embedding $j$.
For sake of simplicity we use the notation
$$\vartheta_\varphi(z):=\vartheta_\varphi(j(z)).$$
With the notation $g=(g_1,g_2)'$,
$h_i=g_{i1}f_1+\cdots+ g_{i8}f_8$ and $\varphi(h)=\phi(g)$ we obtain
for the restriction
$$\vartheta_\varphi(z)=\sum_{h\in\oz(\gz)^2}\varphi(h)
e^{\pii\{ N(h_1)z_1+N(h_2)z_2+\tr(\bar h_1*\Gotz*h_2)\}}.$$
We recall the notion of an
orthogonal modular form of weight $k\in\gz$ with
respect to a subgroup $\Gamma\subset\O^+(\hbox{\rm II}_{2,10})$ of finite index
and with respect to a character $v:\Gamma\to\cz^\bullet$.
It is a  
function $f:\tilde\calH_{10}\to \cz$ with
the properties
$$\eqalign{f(\gamma z)&=v(\gamma)f(z),\cr
f(tz)&=t^{-k}f(z).\cr}$$
The vector space of holomorphic forms
is denoted by $[\Gamma,k,v]$
and by $[\Gamma,k]$ when $v$ is trivial.
These are finite dimensional spaces.
\smallskip
A modular form $f$ is determined by the function
$$ F(z):=F(z_0,z_2,\gotz):=f(1,*,z_0,z_2,\gotz)$$
It satisfies the transformation formula
$$F(\gamma(z_0,z_2,\gotz))=a(\gamma,(z_0,z_2,\gotz))^kF(z_0,z_2,\gotz).$$
Here $\gamma(z_0,z_2,\gotz)$ and $a(\gamma,(z_0,z_2,\gotz))$ are defined as follows:
Consider
$$\gamma((1,*,z_0,z_2,\gotz)=t(1,*,w_0,w_2,\gotw)$$
and define
$$a(\gamma,(z_0,z_2,\gotz))=t^{-1},\quad \gamma(z_0,z_2,\gotz)=(w_0,w_2,\gotw).$$
With these notations we obtain
\proclaim
{Proposition}
{The theta series
$\vartheta_\varphi(z)$ are modular forms of weight 4
for the group $\O^+(\hbox{\rm II})[2]$.}
TsM%
\finishproclaim
Let $f$ be a modular form on some congruence group $\Gamma$
with trivial multiplier system.
We defined in [FS] the value of $f$
at a non-zero rational isotropic vector in $\rz^4\times\oz$  and this value
is $\Gamma$-invariant.
We always can restrict to the case, where
the first component of the isotropic vector is different from 0 because for
any congruence group $\Gamma$ there is a $\Gamma$-equivalent with this property.
Now let
$$R=\pmatrix{r_1&\gotr\cr\bar{\gotr}&r_2},\quad r_1,r_2\in\qz,\ \gotr\in\oz(\qz),$$
be a rational octavic Hermitian matrix. Then we can define the isotropic element
$r=(1,*,r_1,r_2,\gotr)$.  We use the notation $f(R)$ for this value.
\proclaim
{Lemma}
{Let $R$ be a rational octavic Hermitian matrix. The value of
$\vartheta_\varphi(s,Z)$ at  $R$  is
$$\vartheta_\varphi(s,R)=s^{-8}N^{-16}
\sum_{g\mod N}\varphi(g)e^{\pii sR[g]}.$$
Here $N$ denotes a natural number such that $\varphi(g)e^{\pii sR[g]}$
is periodic with respect to $N\oz(\gz)\times N\oz(\gz)$.}
vaC%
\finishproclaim
{\it Proof.\/} Again one can use the modular embedding into the Siegel case
of degree 16.\qed
\smallskip
Recall that we defined the value of a modular form as follows.
One takes for the representative of an isotropic line a primitive vector
in the lattice II$_{2,10}$. This is well-defined up to sign. Hence
we can define the value at a cusp only for forms of even weight.
When we start with an rational Hermitian matrix $R$ the isotropic
vector $r=(1,*,R)$ is usually integral and hence has to be normalized.
One can around the problem of treating this normalizing factor as follows.
Consider the easiest case
$$\vartheta(z):=\vartheta_1(1,2z)=\sum_{h\in\oz(\gz)^2}
e^{2\pii\{ N(h_1)z_1+N(h_2)z_2+\tr(\bar h_1*\Gotz*h_2)\}}.$$
This series has been introduced by Eie and Krieg and they proved that it is
an additive lift.  From this it follows that $\vartheta$ is a modular form
with respect to the full $\O^+(\hbox{\rm II}_{2,10})$. This implies that
the values at all cusps is the same, namely one. Introducing the Gauss sum
$$G(R)=N^{-16}\sum_{g\mod N}e^{\pii R[g]}$$
we obtain now:
\proclaim
{Lemma}
{Let $R$ be a rational octavic Hermitian matrix. The value of
$\vartheta_\varphi(s,Z)$ at the associated cups class is
$${\vartheta_\varphi(s,R)\over G(2R)}={\displaystyle
\sum_{g \mod N}\varphi(g)e^{\pii sR[g]}\over\displaystyle
s^8\sum_{g\mod N}e^{2\pii R[g]}}.$$}
waS%
\finishproclaim
By the way, it is not clear from this formula (but follows from our deduction)
that this expression depends only on the cusp class. This has to do with
rules concerning Gauss sums as reciprocity laws.
\smallskip
We want to compute the values at the cusps in the case $s=1/2$ and
where $\varphi(g)$ only depends on $g$ mod 2.
Theta series of this species sometimes
are called "`of second kind"'.
\proclaim
{Lemma}
{The theta series of second kind
$$\sum_{h\in\oz(\gz)^2}\varphi(h)
e^{\pii\{ N(h_1)z_1+N(h_2)z_2+\tr(\bar h_1*\Gotz*h_2)\}},\quad
\varphi:(\oz(\gz)/2\oz(\gz))^2\lo\cz,$$
are modular forms on the congruence group of level two
$\O^+(\hbox{\rm II}_{2,10})[2]$ (with trivial multipliers).}
FSz%
\finishproclaim
{\it Proof.\/} We recall the construction of the two-fold covering
of the orthogonal group. One has to consider the Clifford algebra
(over $\gz$) of the lattice II$_{2,10}$. The spin group
$\Spin(\hbox{\rm II}_{2,10})$ is a subgroup of its unit group.
There is a natural homomorphism
$$\Spin(\hbox{\rm II}_{2,10})\lo\SO^+(\hbox{\rm II}_{2,10}).$$
It is known that this homomorphism
is surjective. The modular embedding
which we mentioned already should be understood as a homomorphism
$$\Spin(\hbox{\rm II}_{2,10})\lo \Sp(16,\gz).$$
One can reduce this homomorphism mod 2. The spin group and the group
$\O^+$ agree over the filed of two elements. Hence we get a commutative
diagram
$$\matrix{\Spin(\hbox{\rm II}_{2,10})&\lo &\Sp(16,\gz)\cr
\downarrow&&\downarrow\cr
\O^+(\fz_2^{12})&\lo&\Sp(16,\fz_2)}.$$
This shows that a two fold covering of $\O^+(\hbox{\rm II}_{2,10})[2]$
appears as subgroup of a Siegel modular group. This reduces \FSz\
to a Siegel analogue case which is assumed to be known.\qed
\smallskip
We want to compute the values of the theta series of second kind explicitely.
Recall that the cups in the level two case are given by non zero isotropic vectors
from $\fz_2^{12}$. We can take a representative in II$_{2,10}$ such that the first
coordinate is one or two. In terms of the octavic Hermitian matrices $R$
this means that it is sufficient to assume that $R$ either contains only entries
$0$ and $1$ or $0$ and $1/2$.
We prefer to write $R/2$ instead of $R$ with an integral $R$.
Then we have:\smallni
{\it Let $R$ be an octavic Hermitian matrix, such that $R$ either contains $0$ and $1$
or $0$ and $2$.
The value of the second
kind theta function $\sum_{g\;\hbox{\sevenrm integral}}\varphi(g)e^{\pii Z[g]}$
at the cusp class defined by $R/2$ is
$${\displaystyle\sum_{g \mod 4}\varphi(g)e^{\pii R[g]/2}\over\displaystyle
\sum_{g\mod 2}e^{\pii R[g]}}.$$
}
We are going to compute
$$G(\varphi,R)=\sum_{g \mod 4}\varphi(g)e^{\pii R[g]/2}$$
in the cases where $R$ either contains only 0 and 2 or 0 and 1.
This sum can be computed easily:
$$\eqalign{
G(\varphi,R)&=\sum_{g\mod 2}\varphi(g)\sum_{h\mod 2}e^{\pii R[g+2h]/2}\cr&=
\sum_{g\mod 2}\varphi(g)e^{\pii R[g]/2}\sum_{h\mod 2}e^{\pii \tr(\bar h'Rg)}.\cr}$$
The last sum is a sum a sum over the values of a character, which is not zero if and
only if the character is trivial. This means that $Rg$ is a even (contained in
$2(\oz(\gz)\times\oz(\gz))$). The result is
$$\sum_{g\mod 2,\ Rg\;\hbox{\sevenrm even}}\varphi(g)e^{\pii R[g]/2}.$$
Using this formula, the values at the cusps can be calculated.
This can be done
for all the characteristic functions $\varphi$.
Since the weight 4 is the singular weight, we can prove now by calculation:
\proclaim
{Proposition}
{The space of modular forms  $\vartheta_\varphi(z)$ has dimension
715.}
DTa%
\finishproclaim
It follows that this space decomposes into the trivial one-dimensional and
the 714-dimensional irreducible representation.
By the uniqueness of $f(z)$ we obtain that the invariant form $\vartheta(z)$ is contained
in the space generated by the $\vartheta_\varphi(z)$.
Obviously
$\vartheta(2z)$ is nothing else but $\vartheta_\varphi(z)$, where $\varphi$
is the characteristic function of the zero element of $(\oz/2\oz)^2$.
It can be seen form the values at the cusps  that $\vartheta(2z)$ is not contained
in the 714-dimensional space.
Thus  the theorem in the introduction now follows.
\proclaim
{Theorem}
{The space of theta functions of second kind $\vartheta_\varphi(z)$
agrees with the $715$-dimensional
additive lift space described in [FS].
It contains the full invariant form $\vartheta(z)$.}
TaA%
\finishproclaim
\neupara{Enriques surfaces }%
Denote by $U$ the unimodular lattice $\gz\times\gz$ with quadratic form
$(x,x)=2x_1x_2$.
Kondo investigated in [Ko] the case of the lattice
$$M=U\oplus \sqrt2 U\oplus (-\sqrt 2 E_8).$$
This case is related to the
moduli space of marked Enriques surfaces, i.e. Enriques surfaces
with a choice of level 2 structure of the Picard lattice.
\smallskip
The group $\O^+(\hbox{\rm II}_{2,10})[2]$ is related to the  lattice
$$L=\sqrt 2\hbox{II}_{2,10}\cong \sqrt 2 U\oplus \sqrt2 U+\oplus (-\sqrt 2 E_8),$$
cf. [FS].
It can be embedded into Kondo's lattice by means of
$$\sqrt 2 U\lo U,\quad \sqrt 2(x_1,x_2)\loma (x_1,2x_2).$$
Hence, according to [FS],  we have   
\proclaim
{Proposition}
{The embedding of $L=\sqrt 2U\oplus \sqrt2 U\oplus (-\sqrt 2 E_8)\cong
\sqrt 2\hbox{II}_{2,10}$ into
$M=U\oplus \sqrt2 U\oplus (-\sqrt 2 E_8)$ defines an embedding of Kondo's\/
$186$-dimensional space of modular forms of weight four into our\/
$714$-dimensional space.}
KU%
\finishproclaim
Hence as immediate consequence of the previous  section we have
\proclaim
{Theorem}
{Kondo's $186$-dimensional space of modular forms of weight four  is   contained  in the space of theta functions of second kind $\vartheta_\varphi(z)$.}
PU%
\finishproclaim
\vskip2cm\noindent
{\bf References}%
\bigskip
\item{[Bo]} Borcherds, R.:
{\it The moduli space of Enriques surfaces and the fake monster Lie superalgebra.}
Topology . {\bf 35}, 699-710 (1996)
\medskip
\item{[EK]} Eie, M., Krieg, A.:
{\it The Maa\ss\ Space on the Half-Plane of Cayley Numbers of Degree Two,}
Math. Z. {\bf 210}, 113-128 (1992)
\medskip
\item{[FH]} Freitag, E., Hermann, C.F.: {\it Some modular
varieties in low dimension,\/} Advances in Math. {\bf 152}, 203-287 (2000)
\medskip
\item{[FS]}  Freitag,\ E., Salvati-Manni, R.:
{\it Modular forms for the even unimodular lattice of signature (2,10),}
J.~Algebraic Geom. {\bf 16}, 753-791 (2007)
\medskip
\item{[KMY]}  Kawaguchi, S.,  Mukai, S., Yoshikawa, K.: 
{\it Resultants and the  Borcherds $\Phi$-function,\/} arXiv 1308.6454v1
\medskip
\item{[Ko]} Kondo, S.: {\it The moduli space of
Enriques surfaces and Borcherds products,\/}
J. Algebraic Geometry  {\bf 11}, 601-627 (2002)
\medskip
\bye